# Solutions of Sequential Conformable Fractional Differential Equations around an Ordinary Point and Conformable Fractional Hermite Differential Equation

Emrah Ünal[a], Ahmet Gökdoğan[b], Ercan Çelik[c]

[a] Department of Elementary Mathematics Education, Artvin Çoruh University, 08100 Artvin, Turkey
emrah.unal@artvin.edu.tr
[b] Department of Mathematical Engineering, Gümüşhane University, 29100 Gümüşhane, Turkey,
gokdogan@gumushane.edu.tr
[c] Department of Mathematics, Atatürk University, 25400 Erzurum, Turkey,
ecelik@atauni.edu.tr

**Abstract**

In this work, we give the power series solutions around an ordinary point, in the case of variable coefficients, homogeneous sequential linear conformable fractional differential equations of order $2\alpha$. Further, we introduce the conformable fractional Hermite differential equations, conformable fractional Hermite polynomials and basic properties of these polynomials.

**Keywords:** Sequential conformable fractional differential equation, ordinary point, conformable fractional Hermite differential equation, conformable fractional Hermite polynomials

## 1. Introduction

The idea of fractional derivative was raised first by L'Hospital in 1695. Since then, related to the definition of fractional derivatives have been many definitions. The most popular ones of this definitions are Riemann-Liouville and Caputo definitions. For Riemann-Liouville, Caputo and other definitions and the characteristics of these definitions, we refer to reader to [1-3].

Recently, Khalil and friends give a new definition of fractional derivative and fractional integral in [4]. This new definition benefit from a limit form as in usual derivatives. They also proved the product rule, the fractional Rolle theorem and mean value theorem. In [5], Abdeljawad improve this new theory. For instance, definitions of left and right conformable fractional derivatives and fractional integrals of higher order (i.e. of order $\alpha > 1$), Taylor power series representation and Laplace transform of few certain functions, fractional integration by parts formulas, chain rule and Gronwall inequality are provided by him. Legendre conformable fractional equation and Legendre fractional polynomials are studied in [6].

In this work, we analyze the existence of solutions around an ordinary point of conformable fractional differential equation of order $2\alpha$. Then, we give solution of Hermite fractional differential equation. For this solution, we obtain Hermite fractional polynomials with certain special initial conditions. Finally, we introduce the basic properties of Hermite fractional polynomials.

## 2. Conformable Fractional Calculus

**Definition 2.1.** [5] Given a function $f:[a,\infty) \to \mathbb{R}$. Then the left conformable fractional derivative of $f$ order $\alpha$ is defined by

$$(T_\alpha^a f)(x) = \lim_{\varepsilon \to 0} \frac{f(x + \varepsilon(x-a)^{1-\alpha}) - f(x)}{\varepsilon}$$

for all $x > a, \alpha \in (0,1]$. When $a = 0$, it is written as $T_\alpha$. If $(T_\alpha f)(x)$ exists on $(a,b)$ then $(T_\alpha^a f)(a) = \lim_{x \to a^+} (T_\alpha^a f)(x)$.

**Definition 2.2.** [5] Given a function $f:(-\infty, b] \to \mathbb{R}$. Then the right conformable fractional derivative of $f$ order $\alpha$ is defined by

$$(^b_\alpha T f)(x) = -\lim_{\varepsilon \to 0} \frac{f(x + \varepsilon(b-x)^{1-\alpha}) - f(x)}{\varepsilon}$$

for all $x < b, \alpha \in (0,1]$. If $(_\alpha T f)(x)$ exists on $(a,b)$ then $(^b_\alpha T f)(b) = \lim_{x \to b^-} (^b_\alpha T f)(x)$.

**Theorem 2.1.** [4] Let $\alpha \in (0,1]$ and $f, g$ be $\alpha$-differentiable at a point $x > 0$. Then

(1) $\frac{d^\alpha}{dx^\alpha}(af + bg) = a\frac{d^\alpha f}{dx^\alpha} + b\frac{d^\alpha g}{dx^\alpha}$, for all $a, b \in \mathbb{R}$

(2) $\frac{d^\alpha}{dx^\alpha}(x^p) = px^{p-\alpha}$, for all $p \in \mathbb{R}$

(3) $\frac{d^\alpha}{dx^\alpha}(\lambda) = 0$, for all constant functions $f(x) = \lambda$

(4) $\frac{d^\alpha}{dx^\alpha}(fg) = f\frac{d^\alpha}{dx^\alpha}(g) + g\frac{d^\alpha}{dx^\alpha}(f)$

(5) $\frac{d^\alpha}{dx^\alpha}(f/g) = \frac{g\frac{d^\alpha}{dx^\alpha}(f) - f\frac{d^\alpha}{dx^\alpha}(g)}{g^2}$

(6) If, in addition, $f$ is differentiable, then $\frac{d^\alpha}{dx^\alpha}(f(x)) = x^{1-\alpha}\frac{df}{dx}(x)$.

**Theorem 2.2.** [5] Assume $f$ is infinitely $\alpha$-differentiable function, for some $0 < \alpha \leq 1$ at a neighborhood of a point $x_0$. Then $f$ has the fractional power series expansion:

$$f(x) = \sum_{k=0}^{\infty} \frac{\left({}^{(k)}T_\alpha^{x_0} f\right)(x_0)(x-x_0)^{k\alpha}}{\alpha^k k!}, \qquad x_0 < x < x_0 + R^{1/\alpha}, \qquad R > 0.$$

Here, $\left({}^{(k)}T_\alpha^{x_0} f\right)(x_0)$ means the application of the fractional derivative $k$ times.

## 3. Conformable Sequential Fractional Differential Equation and Solutions around an Ordinary Point

The most general sequential linear homogeneous (left) conformable fractional differential equation is

$$^{(n)}T_\alpha^a y + a_{n-1}(x)\,^{(n-1)}T_\alpha^a y + \ldots + a_1(x) T_\alpha^a y + a_0(x) y = 0, \tag{1}$$

where $^{(n)}T_\alpha^a y = T_\alpha^a T_\alpha^a \ldots T_\alpha^a y$, n times.

**Definition 3.1.** Let $\alpha \in (0,1]$, $x_0 \in [a,b]$, $N(x_0)$ be a neighborhood of $x_0$ and $f(x)$ be a real function defined on $[a,b]$. In this case $f(x)$ is said to be $\alpha$-analytic at $x_0$ if $f(x)$ can be expressed as a series of natural powers of $(x - x_0)^\alpha$ for all $x \in N(x_0)$. In other word, $f(x)$ can be expressed as following:

$$\sum_{k=0}^{\infty} c_k (x - x_0)^{k\alpha} \qquad (c_k \in R)$$

This series being definitely convergent for $|x - x_0| < \delta$ ($\delta > 0$). $\delta$ is the radius of convergence of the series.

**Definition 3.2.** Let $\alpha \in (0,1]$, $x_0 \in [a,b]$ and the functions $a_k(x)$ be $\alpha$-analytic at $x_0 \in [a,b]$ for $k = 0,1,2,\ldots,n-1$. In this case, the point $x_0 \in [a,b]$ is said to be an $\alpha$-ordinary point of (1). If a point $x_0 \in [a,b]$ is not $\alpha$-ordinary point, then it is said to be $\alpha$ singular.

**Example 3.1.** a) We consider following the conformable fractional differential equations:

$$x^\alpha T_\alpha y - y = 0,$$

$$x^{2\alpha\,(2)} T_\alpha y - 2 x^\alpha T_\alpha y + x^{2\alpha} y = 0$$

Any point $x = x_0 > 0$ is an $\alpha$-ordinary point for the above equations.

b) Let be

$$(x-1)^\alpha T_\alpha y - y = 0,$$

$$(x-1)^{2\alpha\,(2)} T_\alpha y - 2(x-1)^\alpha T_\alpha y + (x-1)^{2\alpha} y = 0.$$

For these equations, any point $x = x_0 > 1$ is an $\alpha$-ordinary point.

**Theorem 3.1.** Let $\alpha \in (0,1]$, let $c_0, c_1 \in R$ and let $x_0$ be an $\alpha$-ordinary point of the equation

$$T_\alpha^{x_0} T_\alpha^{x_0} y + p(x) T_\alpha^{x_0} y + q(x) y = 0. \tag{2}$$

Then, there exists a solution to the equation (2) as

$$y = \sum_{k=0}^\infty c_k (x - x_0)^{k\alpha} \tag{3}$$

for $x \in (x_0, x_0 + \rho)$ with $\rho = \min\{\delta_1, \delta_2\}$ and initial conditions $c_0 = y(x_0)$, $\alpha c_1 = T_\alpha y(x_0)$ where $\delta_1$ and $\delta_2$ are the radius of convergence of $p(x)$ and $q(x)$, respectively.

**Proof.** Since $x_0$ is an $\alpha$-ordinary point of (2), by definition 3.1 and 3.2 we can write

$$p(x) = \sum_{k=0}^\infty p_k (x - x_0)^{k\alpha} \quad (x \in [x_0, x_0 + \delta_1]; \; \delta_1 > 0) \tag{4}$$

and

$$q(x) = \sum_{k=0}^\infty q_k (x - x_0)^{k\alpha} \quad (x \in [x_0, x_0 + \delta_2]; \; \delta_2 > 0). \tag{5}$$

We seek a solution in form (3) of (2). Substituting (3) and its conformable fractional derivatives in (2), then we obtain

$$\sum_{k=0}^\infty \alpha^2 (k+2)(k+1) c_{k+2} (x - x_0)^{k\alpha}$$

$$+ \left( \sum_{k=0}^\infty p_k (x - x_0)^{k\alpha} \right) \left( \sum_{k=0}^\infty \alpha(k+1) c_{k+1} (x - x_0)^{k\alpha} \right)$$

$$+ \left( \sum_{k=0}^\infty q_k (x - x_0)^{k\alpha} \right) \left( \sum_{k=0}^\infty c_k (x - x_0)^{k\alpha} \right) = 0 \tag{6}$$

We also can write

$$\left( \sum_{k=0}^\infty p_k (x - x_0)^{k\alpha} \right) \left( \sum_{k=0}^\infty \alpha(k+1) c_{k+1} (x - x_0)^{k\alpha} \right)$$

$$= \sum_{k=0}^\infty \left( \sum_{j=0}^k \alpha(j+1) p_{k-j} c_{j+1} \right) (x - x_0)^{k\alpha}$$

$$\tag{7}$$

and

$$\left(\sum_{k=0}^{\infty} q_k(x-x_0)^{k\alpha}\right)\left(\sum_{k=0}^{\infty} c_k(x-x_0)^{k\alpha}\right) = \sum_{k=0}^{\infty}\left(\sum_{j=0}^{k} q_{k-j}c_j\right)(x-x_0)^{k\alpha}.$$

(8)

Hence, if we substitute (7) and (8) in (6), we obtain

$$\sum_{k=0}^{\infty}\left[\alpha^2(k+2)(k+1)c_{k+2} + \sum_{j=0}^{k}\alpha(j+1)p_{k-j}c_{j+1} + \sum_{j=0}^{k} q_{k-j}c_j\right](x-x_0)^{k\alpha} = 0.$$

So, the coefficients $c_k$ must satisfy

$$\alpha^2(k+2)(k+1)c_{k+2} = -\sum_{j=0}^{k}[\alpha(j+1)p_{k-j}c_{j+1}+q_{k-j}c_j]. \qquad (9)$$

We show that if the coefficient $c_k$ are defined by (9), for $k \geq 2$, then the series

$$y = \sum_{k=0}^{\infty} c_k(x-x_0)^{k\alpha}$$

is convergent for $x \in (x_0, x_0 + \rho)$. Let us fix $r$ $(0 < r < \rho)$. Since the series in (4) and (5) are convergent for $x \in [x_0, x_0 + r]$, there is a constant $M > 0$ such that

$$|p_{k-j}| \leq \frac{Mr^{j\alpha}}{r^{k\alpha}} \qquad (k \in N_0; 0 \leq j \leq k) \qquad (10)$$

and

$$|q_{k-j}| \leq \frac{Mr^{j\alpha}}{r^{k\alpha}} \qquad (k \in N_0; 0 \leq j \leq k). \qquad (11)$$

Using (10) and (11) in (9), we obtain

$$\alpha^2(k+2)(k+1)|c_{k+2}| \leq \frac{M}{r^{k\alpha}}\sum_{j=0}^{k}[\alpha(j+1)|c_{j+1}| + |c_j|]r^{j\alpha}$$

$$\leq \frac{M}{r^{k\alpha}}\sum_{j=0}^{k}[\alpha(j+1)|c_{j+1}| + |c_j|]r^{j\alpha} + M|c_{j+1}|r^{\alpha}. \qquad (12)$$

Now, we define

$$C_0 = |c_0|, C_1 = |c_1|$$

and $C_k$ by

$$\alpha^2(k+2)(k+1)C_{k+2} = \frac{M}{r^{k\alpha}}\sum_{j=0}^{k}[\alpha(j+1)C_{j+1} + C_j]r^{j\alpha} + MC_{k+1}r^{\alpha} \qquad (13)$$

for $k \geq 2$.

We can see that an induction yields

$$|c_k| \leq C_k, \quad C_k \geq 0, \quad (k = 0,1,2,\ldots)$$

Now, we analyze for what $x$ the series

$$\sum_{k=0}^{\infty} C_k(x - x_0)^{k\alpha} \tag{14}$$

is convergent.

Using (13), we obtain

$$\alpha^2(k)(k+1)C_{k+1} = \frac{M}{r^{(k-1)\alpha}}\sum_{j=0}^{k-1}[\alpha(j+1)C_{j+1} + C_j]r^{j\alpha} + MC_k r^\alpha \tag{15}$$

$$\alpha^2(k)(k-1)C_k = \frac{M}{r^{(k-2)\alpha}}\sum_{j=0}^{k-2}[\alpha(j+1)C_{j+1} + C_j]r^{j\alpha} + MC_{k-1} r^\alpha. \tag{16}$$

From (15) and (16), we find

$$r^\alpha \alpha^2(k)(k+1)C_{k+1} = \alpha^2(k)(k-1)C_k + \alpha Mkr^\alpha C_k + MC_k r^{2\alpha}.$$

Hence,

$$\frac{C_{k+1}}{C_k} = \frac{\alpha^2(k)(k-1) + \alpha Mkr^\alpha + Mr^{2\alpha}}{r^\alpha \alpha^2(k)(k+1)}$$

is obtained. By the help of the ratio test, we have that

$$\lim_{k\to\infty}\left|\frac{C_{k+1}(x-x_0)^{(k+1)\alpha}}{C_k(x-x_0)^{k\alpha}}\right| = \left(\frac{|x-x_0|}{r}\right)^\alpha < 1$$

Thus, the series (14) converges for $x \in (x_0, x_0 + r)$. This implies that the series (3) converges for $x \in (x_0, x_0 + r)$. Since $r$ was any number satisfying $0 < r < \rho$, the series (3) converges for $x \in (x_0, x_0 + \rho)$.

**Example 3.2.** Find the general solution to the equation

$$^{(2)}T_\alpha y - x^\alpha T_\alpha y - y = 0. \tag{17}$$

For $x_0 = 0$ $\alpha$-ordinary point, we seek a solution in the form (3). Substituting (3) and conformable fractional derivatives of (3) in (17), we have

$$c_2 = \frac{1}{2\alpha^2}c_0$$

and

$$c_{k+2} = \frac{k\alpha + 1}{\alpha^2(k+1)(k+2)}c_k \quad k = 1,2,\ldots$$

Hence,

$$c_2 = \frac{1}{2\alpha^2} c_0 \qquad c_3 = \frac{2\Gamma\left(\frac{1-\alpha}{2\alpha}+2\right)}{\alpha.\Gamma\left(\frac{1}{2\alpha}+1\right)\Gamma(4)} c_1$$

$$c_4 = \frac{\Gamma\left(\frac{1}{2\alpha}+2\right)}{2^{-1}\alpha^3\Gamma\left(\frac{1}{2\alpha}+1\right)\Gamma(5)} c_0 \qquad c_5 = \frac{2^2\Gamma\left(\frac{1-\alpha}{2\alpha}+3\right)}{\alpha^2\Gamma\left(\frac{1-\alpha}{2\alpha}+1\right)\Gamma(6)} c_1$$

$$c_6 = \frac{\Gamma\left(\frac{1}{2\alpha}+3\right)}{2^{-2}\alpha^4\Gamma\left(\frac{1}{2\alpha}+1\right)\Gamma(7)} c_0 \qquad c_7 = \frac{2^3\Gamma\left(\frac{1-\alpha}{2\alpha}+4\right)}{\alpha^3\Gamma\left(\frac{1-\alpha}{2\alpha}+1\right)\Gamma(8)} c_1$$

$$\vdots$$

$$c_{2k} = \frac{\Gamma\left(\frac{1}{2\alpha}+k\right)}{2^{1-k}\alpha^{k+1}\Gamma\left(\frac{1}{2\alpha}+1\right)\Gamma(2k+1)} c_0 \qquad c_{2k+1} = \frac{2^k\Gamma\left(\frac{1-\alpha}{2\alpha}+k+1\right)}{\alpha^k\Gamma\left(\frac{1-\alpha}{2\alpha}+1\right)\Gamma(2k+2)} c_1$$

is obtained. The general solution of (17) is found as

$$y(x) = c_0 \sum_{k=0}^{\infty} \left[ \frac{\Gamma\left(\frac{1}{2\alpha}+k\right)}{2^{1-k}\alpha^{k+1}\Gamma\left(\frac{1}{2\alpha}+1\right)\Gamma(2k+1)} \right] x^{2k\alpha}$$

$$+ c_1 \sum_{k=0}^{\infty} \left[ \frac{2^k\Gamma\left(\frac{1-\alpha}{2\alpha}+k+1\right)}{\alpha^k\Gamma\left(\frac{1-\alpha}{2\alpha}+1\right)\Gamma(2k+2)} \right] x^{(2k+1)\alpha}$$

## 4. Conformable Sequential Fractional Hermite Differential Equation and Conformable Fractional Hermite Polynomials

Consider the conformable fractional Hermite differential equation

$$^{(2)}T_\alpha y - 2\alpha x^\alpha T_\alpha y + 2\alpha^2 m y = 0 \qquad (18)$$

where $\alpha \in (0,1]$, $m$ is a real number. If $\alpha = 1$, then equation (18) becomes the classical Hermite differential equation. Let $m$ be a nonnegative integer. $x = 0$ is an ordinary point of (18). Now we seek a solution as in (3) of (18). Substituting (3) and its conformable fractional derivatives in (18), we have

$$c_2 = -m c_0$$

and

$$c_{k+2} = \frac{2(k-m)}{(k+1)(k+2)} c_k \qquad k = 1,2,\ldots$$

Hence,

$$c_2 = (-1)\frac{2m}{2!}c_0 \qquad c_3 = (-1)\frac{2(m-1)}{3!}c_1$$

$$c_4 = (-1)^2 \frac{2^2 m(m-2)}{4!}c_0 \qquad c_5 = (-1)^2 \frac{2^2(m-1)(m-3)}{5!}c_1$$

$$c_6 = (-1)^3 \frac{2^3 m(m-2)(m-4)}{6!}c_0 \qquad c_7 = (-1)^3 \frac{2^3(m-1)(m-3)(m-5)}{7!}c_1$$

$$\vdots$$

$$c_{2k} = (-1)^k \frac{2^k m(m-2)\ldots(m-2k+2)}{(2n)!}c_0 \qquad c_{2k+1} = (-1)^k \frac{2^k(m-1)(m-3)\ldots(m-2k+1)}{(2n+1)!}c_1$$

is obtained. The general solution of (18) is found as

$$y(x) = c_0 + c_0 \sum_{k=1}^{\infty}\left[(-1)^k \frac{2^k m(m-2)\ldots(m-2k+2)}{(2n)!}\right]x^{2k\alpha} + c_1 x^{\alpha}$$

$$+ c_1 \sum_{k=1}^{\infty}\left[(-1)^k \frac{2^k(m-1)(m-3)\ldots(m-2k+1)}{(2k+1)!}\right]x^{(2k+1)\alpha}$$

Now, we pick initial conditions the following as

$$y(0) = c_0 = (-2)^{\frac{m}{2}}(m-1)!!$$

$$T_\alpha y(0) = \alpha c_1 = 0 \to c_1 = 0$$

where

$$(m-1)!! = \begin{cases}(m-1)(m-3)\ldots 3.1 & (m-1) \text{ odd} \\ (m-1)(m-3)\ldots 4.2 & (m-1) \text{ even}\end{cases}$$

If $c_1 = 0$, then all $c_k = 0$ when $k$ is odd. For these initial conditions, the solution is

$$y(x) = (-2)^{\frac{m}{2}}(m-1)!!\left[1 + \sum_{k=1}^{\infty}\left[(-1)^k \frac{2^k m(m-2)\ldots(m-2k+2)}{(2k)!}\right]x^{2k\alpha}\right]. \qquad (19)$$

Specially, for $m = 6$, the solution is

$$y(x) = 64x^{6\alpha} - 480x^{4\alpha} + 720x^{2\alpha} - 120.$$

This is the 6$^{\text{th}}$ order conformable fractional Hermite polynomial. That is

$$H_6^\alpha(x) = 64x^{6\alpha} - 480x^{4\alpha} + 720x^{2\alpha} - 120.$$

To get the odd order Hermite polynomials, we specify the initial conditions

$$y(0) = c_0 = 0,$$

$$T_\alpha y(0) = \alpha c_1 = -\alpha(-2)^{\frac{m+1}{2}}(m)!! \to c_1 = -(-2)^{\frac{m+1}{2}}(m)!!.$$

If $c_0 = 0$, then all $c_k = 0$ when $k$ is even. For these initial conditions, the solution is

$$y(x) = -(-2)^{\frac{m+1}{2}} (m)!! \left[ x^\alpha + \sum_{k=1}^{\infty} \left[ (-1)^k \frac{2^k (m-1)(m-3)\ldots(m-2k+1)}{(2k+1)!} \right] x^{(2k+1)\alpha} \right] \quad (20)$$

Specially, for $m = 5$, the solution is

$$y(x) = 32x^{5\alpha} - 160x^{3\alpha} + 120x^\alpha.$$

This is the 5$^{\text{th}}$ order conformable fractional Hermite polynomial. That is

$$H_5^\alpha(x) = 32x^{5\alpha} - 160x^{3\alpha} + 120x^\alpha.$$

**Property 4.1.**

(I)    $H_m^\alpha(x) = H_m(x^\alpha)$

(II)    $T_\alpha(H_m^\alpha(x)) = 2m\alpha H_{m-1}^\alpha(x)$

(III)    $^{(m)}T_\alpha(H_m^\alpha(x)) = 2^m m! \alpha^m$

(IV)    $H_{m+1}^\alpha(x) = 2x^\alpha H_m^\alpha(x) - 2m H_{m-1}^\alpha(x)$

(V)    $H_{m+1}^\alpha(x) = 2x^\alpha H_m^\alpha(x) - \alpha^{-1} T_\alpha(H_m^\alpha(x))$

(VI)    $H_m^\alpha(x) = (-\alpha^{-1})^m e^{x^{2\alpha}} {}^{(m)}T_\alpha(e^{-x^{2\alpha}})$

(VII)    $\int_{-\infty}^{\infty} H_m^\alpha(x) H_n^\alpha(x) e^{-x^{2\alpha}} d_\alpha(x) = 0 \;\; m \neq n$ and $\alpha = \frac{1}{2j+1} \;\; j \in N$

(VIII)    $\int_{-\infty}^{\infty} H_m^\alpha(x) H_n^\alpha(x) e^{-x^{2\alpha}} d_\alpha(x) = \frac{1}{\alpha} 2^n n! \sqrt{\pi} \;\; m = n$ and $\alpha = \frac{1}{2j+1} \;\; j \in N$

**Proof.**

**(I)** Proof is obvious.

**(II)** Let $m$ be even. Then, $H_m^\alpha(x)$ is found by the help of (19). Applying conformable derivative to (19), we get

$$T_\alpha(H_m^\alpha(x)) = 2m\alpha \left[ -(-2)^{\frac{m}{2}}(m-1)!! \left[ x^\alpha + \sum_{k=1}^{\infty} \left( (-1)^k \frac{2^k(m-2)\ldots(m-2k)}{(2k+1)!} \right) x^{(2k+1)\alpha} \right] \right].$$

For $m - 1$ is odd, we can write

$$T_\alpha(H_m^\alpha(x)) = 2m\alpha H_{m-1}^\alpha(x).$$

Conversely, let $m$ be odd. Then $H_m^\alpha(x)$ is found by the help of (20). The conformable derivative of order $\alpha$ is

$$T_\alpha(H_m^\alpha(x)) = 2m\alpha \left[ (-2)^{\frac{m-1}{2}}(m-2)!! \left[ 1 + \sum_{k=1}^{\infty} \left( (-1)^k \frac{2^k(m-1)(m-3)\ldots(m-2k+1)}{(2k)!} \right) x^{(2k)\alpha} \right] \right].$$

For $m-1$ is even, similarly, we can write

$$T_\alpha(H_m^\alpha(x)) = 2m\alpha H_{m-1}^\alpha(x). \tag{21}$$

**(III)** We use induction to prove. For $m = 1$, property is provided. That is

$$T_\alpha(H_1^\alpha(x)) = 2\alpha.$$

Assume that the property is true for $m = n$. That is

$$^{(n)}T_\alpha(H_n^\alpha(x)) = 2^n n!\, \alpha^n. \tag{22}$$

From (II), we have

$$T_\alpha(H_{n+1}^\alpha(x)) = 2(n+1)\alpha H_n^\alpha(x). \tag{23}$$

Applying conformable fractional derivative of order $\alpha$, $n$ times, we get

$$^{(n+1)}T_\alpha(H_{n+1}^\alpha(x)) = 2(n+1)\alpha\, ^{(n)}T_\alpha(H_n^\alpha(x)). \tag{24}$$

Substituting (22) in (24), we obtain

$$^{(n+1)}T_\alpha(H_{n+1}^\alpha(x)) = 2^{n+1}(n+1)\alpha^{n+1}. \tag{25}$$

Hence, proof is completed.

**(IV)** For $H_{m+1}^\alpha(x)$ is a solution of (18), we can write

$$^{(2)}T_\alpha(H_{m+1}^\alpha(x)) - 2\alpha x^\alpha T_\alpha(H_{m+1}^\alpha(x)) + 2\alpha^2 m H_{m+1}^\alpha(x) = 0. \tag{26}$$

Using (II), we have

$$H_{m+1}^\alpha(x) = 2x^\alpha H_m^\alpha(x) - 2m H_{m-1}^\alpha(x). \tag{27}$$

**(V)** Using (21) and (27), the result is found.

**(VI)** We prove by induction. For $m = 1$, the property is provided.

Assume the property is true for $m = n$. That is

$$H_n^\alpha(x) = (-\alpha^{-1})^n e^{x^{2\alpha}}\, ^{(n)}T_\alpha(e^{-x^{2\alpha}}). \tag{28}$$

If (28) is substituted in the property (V), then

$$H_{n+1}^\alpha(x) = (-\alpha^{-1})^{n+1} e^{x^{2\alpha}}\, ^{(n+1)}T_\alpha(e^{-x^{2\alpha}})$$

is obtained. Hence, proof is completed.

**(VII)** The definition 2.1 is given for $x > 0$. To avoid the problem of being undefined on $(-\infty, 0]$, we assume $\alpha = \frac{1}{2j+1}$, with $j$ any natural number. Since $H_m^\alpha(x)$ and $H_n^\alpha(x)$ are solution for the equation (18), respectively, then

$$^{(2)}T_\alpha\left(H_m^\alpha(x)\right) - 2\alpha x^\alpha T_\alpha\left(H_m^\alpha(x)\right) + 2\alpha^2 m H_m^\alpha(x) = 0 \tag{29}$$

$$^{(2)}T_\alpha\left(H_m^\alpha(x)\right) - 2\alpha x^\alpha T_\alpha\left(H_m^\alpha(x)\right) + 2\alpha^2 m H_m^\alpha(x) = 0 \tag{30}$$

Multiply (29) by $e^{-x^{2\alpha}} H_n^\alpha(x)$ and (30) by $e^{-x^{2\alpha}} H_m^\alpha(x)$ and subtract the resulting equation to get

$$T_\alpha\left[e^{-x^{2\alpha}} T_\alpha\left(H_m^\alpha(x)\right)\right] H_n^\alpha(x) - T_\alpha\left[e^{-x^{2\alpha}} T_\alpha\left(H_n^\alpha(x)\right)\right] H_m^\alpha(x) + 2\alpha^2(m-n) e^{-x^{2\alpha}} H_n^\alpha(x) H_m^\alpha(x) = 0. \tag{31}$$

If we apply the fractional integral to equation (31), then we get

$$\int_{-\infty}^{\infty} \left(T_\alpha\left[e^{-x^{2\alpha}} T_\alpha\left(H_m^\alpha(x)\right)\right] H_n^\alpha(x) - T_\alpha\left[e^{-x^{2\alpha}} T_\alpha\left(H_n^\alpha(x)\right)\right] H_m^\alpha(x)\right) d_\alpha x$$

$$+ 2\alpha^2(m-n) \int_{-\infty}^{\infty} e^{-x^{2\alpha}} H_n^\alpha(x) H_m^\alpha(x) d_\alpha x = 0$$

If we apply fractional integration by parts to the first fractional integral in the above equation, we find the result of this fractional integral as zero. Hence proof is completed.

**(VIII)** The definition 2.1 is given for $x > 0$. To avoid the problem of being undefined on $(-\infty, 0]$, we assume $\alpha = \frac{1}{2j+1}$, with $j$ any natural number. Let us define

$$I_{m,n} = \int_{-\infty}^{\infty} H_m^\alpha(x) H_n^\alpha(x) e^{-x^{2\alpha}} d_\alpha(x)$$

Then

$$I_{n-1,n+1} = \int_{-\infty}^{\infty} H_{n-1}^\alpha(x) H_{n+1}^\alpha(x) e^{-x^{2\alpha}} d_\alpha(x) = 0$$

Using the property (IV), we get

$$\int_{-\infty}^{\infty} H_{n-1}^\alpha(x) \left(2x^\alpha H_n^\alpha(x) - 2n H_{n-1}^\alpha(x)\right) e^{-x^{2\alpha}} d_\alpha(x) = 0.$$

i.e.

$$\int_{-\infty}^{\infty} 2x^\alpha H_{n-1}^\alpha(x) H_n^\alpha(x) e^{-x^{2\alpha}} d_\alpha(x) = 2n I_{n-1,n-1}. \tag{32}$$

Recall that

$$H_n^\alpha(x) = (-\alpha^{-1})^n e^{x^{2\alpha}} {}^{(n)}T_\alpha(e^{-x^{2\alpha}}).$$

Thus, equation (32) becomes

$$-(\alpha^{-1})^{2n-1} \int_{-\infty}^{\infty} 2x^\alpha e^{x^{2\alpha}} \left({}^{(n-1)}T_\alpha(e^{-x^{2\alpha}})\right)\left({}^{(n)}T_\alpha(e^{-x^{2\alpha}})\right) d_\alpha(x) = 2nI_{n-1,n-1}. \quad (33)$$

Note that

$$2x^\alpha e^{x^{2\alpha}} {}^{(n-1)}T_\alpha(e^{-x^{2\alpha}}) = \frac{1}{\alpha} T_\alpha\left[e^{x^{2\alpha}} {}^{(n-1)}T_\alpha(e^{-x^{2\alpha}})\right] - \frac{1}{\alpha} e^{x^{2\alpha}} {}^{(n)}T_\alpha(e^{-x^{2\alpha}}).$$

Then, equation (33) becomes

$$\int_{-\infty}^{\infty} (\alpha^{-1})^{2n} e^{x^{2\alpha}} {}^{(n)}T_\alpha(e^{-x^{2\alpha}}) {}^{(n)}T_\alpha(e^{-x^{2\alpha}}) d_\alpha(x) - (\alpha^{-1})^{2n} \int_{-\infty}^{\infty} T_\alpha\left[e^{x^{2\alpha}} {}^{(n-1)}T_\alpha(e^{-x^{2\alpha}})\right] {}^{(n)}T_\alpha(e^{-x^{2\alpha}}) d_\alpha(x)$$

$$= 2nI_{n-1,n-1}$$

Using the property (IV) and integration by parts for conformable fractional derivative, we have

$$I_{n,n} - (\alpha^{-1})^{2n} \left[e^{x^{2\alpha}} {}^{(n-1)}T_\alpha(e^{-x^{2\alpha}}) {}^{(n)}T_\alpha(e^{-x^{2\alpha}})\right]_{-\infty}^{\infty} +$$

$$\int_{-\infty}^{\infty} (\alpha^{-1})^{2n} e^{x^{2\alpha}} {}^{(n-1)}T_\alpha(e^{-x^{2\alpha}}) {}^{(n+1)}T_\alpha(e^{-x^{2\alpha}}) d_\alpha(x) = 2nI_{n-1,n-1}.$$

We also have following equations:

$$\int_{-\infty}^{\infty} (\alpha^{-1})^{2n} e^{x^{2\alpha}} {}^{(n-1)}T_\alpha(e^{-x^{2\alpha}}) {}^{(n+1)}T_\alpha(e^{-x^{2\alpha}}) d_\alpha(x) = I_{n-1,n+1} = 0,$$

$$-(\alpha^{-1})^{2n} \left[e^{x^{2\alpha}} {}^{(n-1)}T_\alpha(e^{-x^{2\alpha}}) {}^{(n)}T_\alpha(e^{-x^{2\alpha}})\right]_{-\infty}^{\infty} = 0.$$

Hence, recurrence equation

$$I_{n,n} = 2nI_{n-1,n-1}$$

is obtained. Repeating this operation $n$ times yields the result

$$I_{n,n} = 2^n n! I_{0,0}$$

where

$$I_{0,0} = \int_{-\infty}^{\infty} e^{-x^{2\alpha}} d_\alpha(x) = \frac{\sqrt{\pi}}{\alpha}$$

Hence proof is completed.

**5. Conclusion**

In this work, we give power series solutions around an ordinary point in homogenous case of sequential linear differential equation of conformable fractional of order $2\alpha$ with variable coefficients. In addition, solving Hermite fractional differential equation, we obtain Hermite fractional polynomials for certain initial conditions. It is appeared that the results obtained in this work correspond to the results which are obtained in ordinary case.